\documentclass[10pt,a4paper]{amsart}
\usepackage{amssymb}
\usepackage[all]{xy}
\CompileMatrices

\def \R{\mathbb R}

\def \Z{\mathbb Z}
\def \Q{\mathbb Q}

\theoremstyle{plain}
\newtheorem{theorem*}{Theorem}

\newtheorem{lemma*}{Lemma}
\newtheorem{corollary*}[theorem*]{Corollary}

\begin{document}
\begin{center}
\bf{ON RECIPROCITY}
\end{center}

\vskip 0.4cm

\begin{center}
{\footnotesize{FLORIAN DELOUP{\footnote{Supported by a E.U.
Research Grant MERG-CT-2004-510590.}} AND VLADIMIR TURAEV}}
\end{center}

\vskip 0.3cm

\noindent {\footnotesize{Abstract. We prove a reciprocity formula
between Gauss sums that is used in the computation of certain
quantum invariants of 3-manifolds. Our proof uses the discriminant
construction applied to the tensor product of lattices.}}

\vskip 0.7cm

%
%

%
%
%
%
%
%

%
%
%

%
%
%

In 1992, L. Jeffrey \cite{Jeffrey} (\S 4.2, Proposition 4.3, p.
583) stated the following reciprocity formula
 for Gauss sums. Let $W$ be  a lattice  of finite rank $l$ and $\langle \cdot, \cdot
\rangle$ be an inner product in the real vector space $W_{\R}=W
\otimes_{\Z} \R$. Let  $W^{\bullet}  = \{ x \in W_{\R} \ | \
\langle x, W \rangle \subset  \Z \}$ be the dual lattice.  Let
$h:W_{\R} \to W_{\R}$ be a self-adjoint automorphism, $\psi \in
W_{\R}$ and $r $ a positive integer. Assume that $h(W^{\bullet})
\subset  W^{\bullet}$ and
\begin{equation} \label{eq:hypo} \begin{split}
\frac{r}{2} \langle x, h(x) \rangle,\    \langle x, h(x') \rangle,
\  r\langle x, \psi \rangle \in \Z, \ \ {\hbox{for all}}\ x, x'
\in W, \\
\frac{r}{2} \langle y, h(y) \rangle,\    \langle y, ry' \rangle,\
 r \langle y, \psi \rangle \in \Z, \ \ {\hbox{for all}}\ y, y'
\in W^{\bullet}. \end{split}
\end{equation}
The reciprocity formula     states that
\begin{equation} \label{eq:rec-ht}
{\rm{vol}}(W^{\bullet}) \sum_{x \in W/rW}   \exp \Bigl( \frac{\pi
i}{r}(\langle \tilde x, h(\tilde  x) \rangle
  + 2 r \langle \tilde{x}, \psi \rangle) \Bigr)  \end{equation}
$$ = \biggl( \det \left(\frac{h}{i}\right) \biggr)^{-1/2} r^{l/2}
\sum_{y \in W^{\bullet}/h(W^{\bullet})} \exp\left(-\pi i r \langle
\tilde  y + \psi, h^{-1}(\tilde  y + \psi) \rangle \right)
$$
where $\tilde  x \in W$ stands for an arbitrary lift of $x$ and
$\tilde y\in W^{\bullet}$ stands for an arbitrary lift of $y$.
(The sums above do not depend on the choice of the lifts.) Note
that   the volume vol$(L) $  of a lattice $L\subset W_{\R}$ is the
absolute value of the determinant of a matrix expanding  a basis
of $L$ with respect to an orthonormal basis of $W_{\R}$.

Formula (\ref{eq:rec-ht}) was later used by several authors, see
for instance \cite[p.\ 623]{HansenTakata} and \cite[Th. 8.1,
formulas (8.4), (8.5)]{MPR}. However, Formula (\ref{eq:rec-ht})
has an
 obvious sign ambiguity  in the choice of a square root of $\det
 \left( {h}/{i}\right) $. Moreover,  one important factor is
  missing as is clear from the following corrected version of
  (\ref{eq:rec-ht}).

Define a symmetric bilinear form $g:W \times W \to \Z$ by
 $ g(x,y) = \langle x, h(y) \rangle$ for $ x, y \in W$.
Recall that the  signature $\sigma(g) \in \Z$ is defined as the
number of positive terms  minus the number of negative terms in a
diagonal matrix presenting the bilinear extension $ W_\R \times W_\R
\to \R$ of $g$. Our main result is

\begin{theorem*}
\begin{equation} \label{2eq:rec-ht}  {\rm{vol}}(W^{\bullet}) \sum_{x \in W/rW}   \exp
 \frac{\pi i}{r}\Bigl(\langle \tilde x, h(\tilde  x)\rangle
  + 2 r \langle {\tilde{x}}, \psi \rangle \Bigr)   \end{equation}
$$= |\det(h)|^{-1/2} \exp\left(  \pi i\frac{\sigma(g)}{4}\right)\
r^{l/2} \sum_{y \in W^{\bullet}/h(W^{\bullet})} \exp\left(-\pi i r
\langle \tilde  y + \psi, h^{-1}(\tilde  y + \psi) \rangle \right)
$$
where $\tilde  x \in W$, $\tilde y\in W^{\bullet}$ are arbitrary
lifts of $x,y$ as above.
\end{theorem*}

This formula belongs to a long sequence of reciprocity formulas
for Gauss sums due to Cauchy, Dirichlet, Kronecker, Krazer and
others. In 1997, in his Ph.\ D.\ thesis \cite{De0} (see also
\cite{De1}), the first author established a reciprocity formula
for Gauss sums generalizing Krazer's formula. This was further
generalized by the second author  in \cite{Tu}. The aim of this
note is to  deduce Formula (\ref{2eq:rec-ht}) from our reciprocity
formulas.

 In Section 1 we recall our
reciprocity formulas and deduce them from a classical van der Blij
formula. In Section 2  we prove  Formula (\ref{2eq:rec-ht}).

\section{General reciprocity formulas}

\subsection{The Gauss sums}
A {\emph{quadratic function}} on a finite abelian group $G$ is a map
$q:G \to \Q/\Z$ such that $b_q(x,y) = q(x+y) - q(x) - q(y)$ is
$\Z$-bilinear in $x$ and $y$.  For a group $H\subset G$, set
$H^\bot= \{x\in G\, \vert\, b_q (x,H)=0\}$ and $ |H|=$ card $(H)$.
We say that  $q$ is {\emph{nondegenerate}} if $G^{\bot} = 0$. Let
$$\gamma(G,q) = |G|^{-1/2} |G^{\bot}|^{-1/2} \sum_{x \in G} e^{2
\pi i q(x)}$$ be the Gauss sum associated to $q$  and normalized
so that its absolute value    is either $0$ or $1$. The latter
case occurs if and only if $q|_{G^{\bot}}= 0$.

We begin with two classical lemmas.  In both lemmas, $q:G \to
\Q/\Z$ is a nondegenerate quadratic function. For the first lemma,
see for example \cite{Scharlau},  pp.168 and 171. The proofs are
given there for homogeneous $q$ (in the sense that $q(nx)=n^2
q(x)$ for any  $n\in  \Z, x\in G$), but are valid for all  $q$.

\begin{lemma*} \label{lem:classical}
For any  group $H \subset  G$   such that $H \subset H^{\bot}$ and
$q\vert_H=0$, the function $q$ induces a quadratic function
$q_H:H^{\bot}/H \to \Q/\Z$.  Moreover, $\gamma(G,q) =
\gamma(H^{\bot}/H, q_H)$.
\end{lemma*}

In the following lemma (\cite{Springer}, \S 1.7), bar denotes
complex conjugation.

\begin{lemma*} \label{th:rec-gen}  For any group $A\subset G$,
\begin{equation}
\gamma(A, q|_{A}) = \gamma(G, q) \, \overline{\gamma(A^\bot,
q|_{A^\bot})}. \label{eq:rec-gen}
\end{equation}
\end{lemma*}

 \noindent{\bf{Proof}}. Set $B=A^\bot$ and $H = A \cap B$. If
$q|_{H} \not= 0$, then $\gamma(A, q|_{A}) = \gamma(B, q|_{B}) = 0$
and Formula (\ref{eq:rec-gen}) holds. Suppose that  $q|_{H} = 0$.
Since $q$ is non-degenerate, $|G| / | H^{\bot}|= |G/H^{\bot} |= |$
Hom $(H, \Q/\Z) |= |  H |$ and similarly $ |G|/ | B|= |A |$.
Therefore $|A+B| =|A| |B|/ |H| =  |G|/ | H|= |H^{\bot} |$.  The
obvious inclusion $A+B\subset H^{\bot} $ implies that $A+B
=H^{\bot}$. By Lemma \ref{lem:classical}, $q$ induces a quadratic
function $q_H$ on $H^{\bot}/H$ and $\gamma(G,q) =
\gamma(H^{\bot}/H, q_H)  $. Observe that  $(H^{\bot}/H, q)$ is an
orthogonal sum of the induced quadratic functions on $A/H$ and
 $B/H$. Their   Gauss sums are equal to    $\gamma(A, q|_{A})$ and
 $\gamma(B, q|_{B})$ respectively. By     multiplicativity with
respect to the orthogonal  sum,
$$ \gamma(G,q) =
\gamma(H^{\bot}/H, q)  = \gamma(A, q|_{A}) \, \gamma(B, q|_{B}).$$
Formula (\ref{eq:rec-gen}) follows. \hfill $\blacksquare$

\subsection{The discriminant}\label{subsec:lattices} To compute   the Gauss sums explicitly, it is
convenient to use the   discriminant construction on lattices
which we now recall. A {\emph{lattice}}   is a finitely generated
free abelian group. A \emph{ bilinear lattice} $(V,f)$ is a
symmetric nondegenerate bilinear form $f: V\times V \to \Z$ on a
lattice $V$.  A  bilinear lattice $(V,f)$ extends  to a  symmetric
nondegenerate bilinear form $f_{\Q}:V_{\Q} \times V_{\Q} \to \Q$
where  $V_{\Q} = V \otimes \Q$ is  a vector space  over $\Q$ of
finite dimension equal to the rank of $V$.   (Here and below
$\otimes =\otimes_{\Z}$). Let $V^\sharp=\{ x\in V_{\Q} \ | \
f_{\Q}(x,V) \subset  \Z\}$ be the   dual lattice.   Clearly $V
\subset  V^{\sharp}$. Consider the finite abelian group $G_f =
V^{\sharp}/V$ and define a symmetric bilinear form $\lambda_f: G_f
\times G_f \to \Q/\Z$  by
\begin{equation}
\label{eq:pairing} \lambda_f\left(x (\hbox{mod}\ V) , y
(\hbox{mod}\ V) \right) = f_{\Q}(x,y)\ \hbox{mod}\ \Z,
\end{equation}
where $x, y \in V^{\sharp} $. Note that $\lambda_f$ is
nondegenerate.

A \emph{Wu class} (resp.\ \emph{integral Wu class})
  for $(V,f)$ is an element $v \in V_{\Q}$ (resp.   $v \in V$) such that
$ f(x,x)-f_{\Q}(x,v) \in 2\Z$ for all $ x\in V$. The following
lemma is well-known (see for example \cite{LW}, Lemma 1.6 (i)).

\begin{lemma*}
Every bilinear lattice $(V,f)$ has an integral Wu class.
\end{lemma*}

\noindent{\bf{Proof}}. Set $V'=V \otimes \Z/2\Z$ and consider the
induced pairing $f':V' \times V' \to \Z/2\Z$. If $v' \in V'$
verifies $f'(x,x) = f'(x,v')$ for all $x \in V'$, then any lift of
$v'$ to $V$  is an integral Wu class for $f$. The map $x \mapsto
f'(x,x)$ is a homomorphism $V' \to \Z/2\Z$. If $f'$ is
nonsingular, that is  the adjoint map $V' \to
{\rm{Hom}}(V',\Z/2\Z)$ is an isomorphism, then  there obviously
exists $v' \in V'$ as required. The case of singular  $f'$ reduces
to the nonsingular  case by splitting $(V',f')$ as an orthogonal
sum of a nonsingular pairing and a zero
pairing. \hfill $\blacksquare$\\

For a Wu class $v \in V_{\Q}$, define a quadratic function
$\varphi_{f,v}: G_f \to \Q/\Z$   by
\begin{displaymath}
\varphi_{f,v}\left(x (\hbox{mod}\
V)\right)=\frac{1}{2}(f_{\Q}(x,x)-f_{\Q}(x,v))\ (\hbox{mod}\ \Z)
\end{displaymath}
for $x \in  V^{\sharp}$. One calls  $  \varphi_{f,v} $  the
\emph{discriminant} of $(V,f,v)$. Its associated bilinear pairing
is $\lambda_{f}$. The discriminant construction
  preserves orthogonal sums and produces all
nondegenerate quadratic functions on  finite abelian groups. For
more on this, see \cite{LW}, \cite{Durfee},  \cite{DM}. The
following  formula due to van der Blij \cite{bl} computes the
Gauss sum $\gamma(G_{f}, \varphi_{f, v})$ via the  signature
$\sigma(f) \in \Z$.

\begin{lemma*}
\label{lem:vanderblij} $\gamma(G_{f}, \varphi_{f,v}) = e^{ {\frac{
\pi i}{4}}  ( \sigma(f) - f_{\Q}(v,v)  )}$.
\end{lemma*}

\subsection{Reciprocity   for tensor products}
\label{subsec:tensor} A  study of the discriminant of a tensor
product  leads to a  reciprocity formula \cite[\S 1.3]{Tu} which
we now state. Let $(V,f)$ and $(W,g)$ be
 bilinear lattices. Set $Z=V\otimes W$ and define a   (symmetric nondegenerate)
bilinear pairing  $f \otimes g: Z \times Z \to \Z$   by $$(f
\otimes g)(x \otimes y, x' \otimes y') = f(x,x') \, g(y,y')\ \ \
{\rm{for}} \ x, x' \in V, \ y, y' \in W.$$ Let $z \in Z \otimes
\Q$ be a Wu class for
  $(Z, f \otimes
g)$. The discriminant of $(Z, f \otimes g, z)$ is a nondegenerate
quadratic function $\varphi_{f \otimes g, z}: G_{f \otimes g} \to
\Q/\Z$ where $G_{f \otimes g}=Z^{\sharp}/Z$. Define a homomorphism
$j_{f}:G_{f} \otimes W \to G_{f \otimes g}$   by
$$j_{f}\left( x\, (\hbox{mod}\ V)  \otimes y \right) =   x \otimes y \,\,(\hbox{mod}\ Z) $$
where $x\in V^{\sharp}, y\in W$.  Similarly, define a homomorphism
$j_{g}:V \otimes G_{g} \to G_{f \otimes g}$ by
$$j_{g}\left( x \otimes y \, (\hbox{mod}\ W) \right) =  x \otimes y \,\, (\hbox{mod}\ Z) $$
where $x\in V, y\in W^{\sharp}$.

\begin{theorem*} \label{th:rec-tensor}
 $$
\gamma(G_{f} \otimes W, \varphi_{f \otimes g, z} \circ j_{f})
 =  e^{
{\frac{ \pi i}{4}}  ( \sigma(f) \sigma(g)  - (f\otimes
g)_{\Q}(z,z) )}  \, \,   \overline{\gamma(V \otimes G_{g},
\varphi_{f \otimes g,z} \circ j_{g})}.$$
\end{theorem*}

\noindent{\bf{Proof}}.  It is easy to check  that $j_{f}$ and
$j_{g}$ are one-to-one, the groups $A=j_{f}(G_{f} \otimes W)$ and
$B=j_{g}(V \otimes G_{g})$ are orthogonal in $G_{f \otimes g}$,
and  ${\rm{Coker}}\ j_{f} = V^{*} \otimes G_{g}$ where $V^{*}=$
Hom $(V, \Z)$ (cf.\ \cite[\S 3]{Tu}).  Therefore  $|G_{f \otimes
g}| = |A| \ |V^{*} \otimes G_{g}| = |A|\ |B| $. Since $|G_{f
\otimes g}| = |A| \ |A^{\bot}|$, we deduce $B=A^\bot$. Lemma
\ref{lem:vanderblij} and the equality $\sigma(f \otimes g) =
\sigma(f) \,\sigma(g) $ give
$$\gamma(G_{f \otimes g}, \varphi_{f \otimes
g, z}) = e^{ {\frac{ \pi i}{4}}  ( \sigma(f) \sigma(g)  -
(f\otimes g)_{\Q}(z,z) )}.$$  The    desired formula follows then
from
Lemma  \ref{th:rec-gen}. \hfill $\blacksquare$\\

We  describe two special cases of Theorem \ref{th:rec-tensor}
first established in \cite[\S 1.4]{De1} and \cite[Lemma 1.4]{De0}.
Let  $v \in V$ and $w \in W$ be integral Wu classes for $f$ and
$g$, respectively.  Clearly, $ v\otimes w\in V\otimes W$  is a Wu
class for $f\otimes g$.  Recall that for  a quadratic function
$q:G \to \Q/\Z$ there is a unique quadratic function   $ q \otimes
g :G \otimes W \to \Q/\Z$ such that    $(q \otimes g)(x \otimes y)
=q(x) g(y,y)$ for  $x \in G, y \in W$.  The associated bilinear
form of $q \otimes g$ is $b_q \otimes g$. Similarly, one has a
quadratic function $ f \otimes q :V \otimes G \to \Q/\Z$.
  Applying Theorem
\ref{th:rec-tensor} to $z=v\otimes w$ and observing that
$$ \varphi_{f \otimes g, v
\otimes w} \circ j_{f} = \varphi_{f, v} \otimes g\ \ \ {\rm{and}}\
\ \ \varphi_{f \otimes g, v \otimes w} \circ j_{g} = f \otimes
\varphi_{g,w} $$ we obtain
$$
\gamma (G_{f} \otimes W,  \varphi_{f,v} \otimes g )   = e^{ \frac{
\pi i}{4}(\sigma(f) \sigma(g) - f(v,v)g(w,w))
  } \,\,   \overline{\gamma (V \otimes G_{g}, f \otimes
\varphi_{g,w}   )}.
 $$
Each pair  $a\in V^{\sharp}, b\in W^{\sharp}$ determines a
homomorphism $p_1(a\otimes b) : G_{f} \otimes W \to \Q/\Z$
sending $x\otimes y$ to $\lambda_f( a (\hbox{mod}\ V), x)\, g_{\Q}
(b,y)$ for $x\in G_f, y\in W$. Similarly, the homomorphism $p_2
(a\otimes b) : V\otimes G_{g}  \to \Q/\Z$   sends $x\otimes y$ to
$f_{\Q} (a,x)\, \lambda_g( b (\hbox{mod}\ W), y) $ for $x\in V,
y\in G_g$. This extends by linearity to homomorphisms $p_1:
V^{\sharp} \otimes W^{\sharp} \to $ Hom $(G_{f} \otimes W,
\Q/\Z)$ and $p_2: V^{\sharp} \otimes W^{\sharp} \to $ Hom $(V
\otimes G_{g},  \Q/\Z)$.

\begin{theorem*} \label{th:tensor-del}
 For any $\zeta \in V^{\sharp} \otimes W^{\sharp}$,
$$
\gamma (G_{f} \otimes W,  \varphi_{f,v} \otimes g  + p_1(\zeta ) )
$$
$$  = e^{\frac{  \pi i}{4}(\sigma(f) \sigma(g) - f(v,v)g(w,w)  - 4 (f_{\Q} \otimes g_{\Q})(\zeta,
\zeta-v\otimes w))} \,\, \overline{\gamma (V \otimes G_{g}, f
\otimes \varphi_{g,w} + p_2 (\zeta) )}.
 $$
\end{theorem*}

\noindent{\bf{Proof}}.  Apply   Theorem \ref{th:rec-tensor} to the
Wu class $z=v\otimes w  - 2\zeta$ of  $f\otimes g$.
 \hfill $\blacksquare$

\section{Proof of Formula (\ref{2eq:rec-ht})}

\subsection{A  particular case} \label{subsec:further}
We begin with a particular case of Theorem \ref{th:tensor-del}.
Let $V=\Z$ and
 $f:V\times V \to \Z$ be   defined by $f(x,y) = r x y$   for an integer $r\geq 1$.
  Clearly, $ V^{\sharp}=  \frac{1}{r} \Z\subset
\Q$ and $G_f= V^{\sharp}/V$ is a cyclic group of order $r$
generated by $1/r  (\hbox{mod}\ V)$. For any
 $x, y\in V^{\sharp}$ and an integral Wu class $v \in
V$,
$$\lambda_f(x(\hbox{mod}\ V), y(\hbox{mod}\ V))=rxy \
({\rm{mod}}\ \Z),\,\,\,\,\,  \varphi_{f,v}(x (\hbox{mod}\ V)) =
\frac{ x (x-v)   }{2r} \ ({\rm{mod}}\ \Z).$$ In the sequel  $v =
0\in V$ for even $r$ and   $v = 1\in  V$ for odd  $r$.

Let $(W,g) $ be a  bilinear  lattice of rank $l$ with integral Wu
class $w\in W$. We shall assume that in the case of  odd $r$, the
form $g$ is even and $w=0$. Clearly,  $G_{f} \otimes W=W/rW$ and
$\lambda_{f} \otimes  g$ is the composition of the bilinear
pairing $g_r: W/rW \times W/rW \to \Z/r\Z$  induced by $g$ with
the embedding $\Z/r\Z\hookrightarrow \Q/\Z$ sending $1
({\rm{mod}}\, r)$ to $1/r\, ({\rm{mod}}\, \Z)$. The  annihilator
$T=(W/rW)^{\bot}$ of  $\lambda_{f} \otimes g$ coincides with the
annihilator of $g_r$. Fix $\zeta \in V^{\sharp} \otimes
W^{\sharp}=   \frac{1}{r} W^{\sharp} \subset W\otimes \Q$.  It
follows from the definitions and the assumptions on $g$ that the
quadratic function $\varphi_{f,v} \otimes g  + p_1(\zeta )$ on
$W/rW $  sends $x\in W/rW$ to $\frac{1}{2r} g(\tilde x, \tilde x)+
g_{\Q}(\tilde x,\zeta) (\hbox{mod}\ \Z)$ where $\tilde x\in W$ is
any lift of $x$. Since $|W/rW| = r^{l}$,
$$\gamma (G_{f} \otimes W,  \varphi_{f,v} \otimes g  + p_1(\zeta ) )
 =    r^{-l/2}\, |T|^{-1/2} \sum_{x \in W/rW} e^{\frac{\pi i}{r}
(g(\tilde x, \tilde x) +  g_{\Q}(\tilde x,  2 r \zeta))} .
$$

 We have $V \otimes G_{g}=G_g$ and $f\otimes \lambda_g= r \lambda_g$. Since $\lambda_g$ is non-degenerate,
the annihilator $G_{g}^{\bot }$ of $f\otimes \lambda_g$ is equal
to $\{ y \in G_{g} \ | \ ry = 0 \}$. The latter group is
isomorphic to $T$ via $y (\hbox{mod}\ W)\mapsto ry (\hbox{mod}\ r
W)$ for any $y\in W^{\sharp}$ with  $ry\in W$. The   quadratic
function $f \otimes \varphi_{g,w} + p_2 (\zeta) : G_g\to \Q/\Z$
sends $y\in G_g=W^{\sharp}/W$ to $$r \varphi_{g,w}(y,y) +
\lambda_{g}(y, r\zeta (\hbox{mod}\ W))=(r/2) g_{\Q}(\tilde y,
\tilde y -w) + r g_{\Q}(\tilde y,\zeta) (\hbox{mod}\ \Z)$$ where
$\tilde y\in W^{\sharp}$ is any lift of $y$. If $r$ is odd, then
$w = 0$ so $\frac{r}{2}g_{\Q}(\tilde{y}, w)  = 0$. If $r$ is even,
then $ {r}/{2} \in  \Z$ and  $g_{\Q}(\tilde{y},w) \in \Z$. Hence,
in all cases, $({r}/{2}) g_{\Q}(\tilde{y},w)   = 0({\rm{mod}}\
\Z)$. Therefore
$$\gamma (V \otimes G_{g}, f \otimes
\varphi_{g,w} + p_2 (\zeta) )=    |G_g|^{-1/2}\,  |T|^{-1/2}
\sum_{y \in W^{\sharp}/W} e^{\pi i r  g_{\Q}(\tilde y,\tilde y +
2\zeta)}.
$$
Since $\sigma(f) =1$ and $v\otimes w=0$,  Theorem
\ref{th:tensor-del} and the computations above give
\begin{equation} \label{eq:special1}   r^{-l/2}\,   \sum_{x \in W/rW} e^{\frac{\pi i}{r}
 (g_{\Q}(\tilde x, \tilde x) + 2 r \langle \tilde x, \zeta \rangle )} \end{equation}
$$ =e^{\frac{  \pi i}{4}\sigma(g) }{|G_g|^{-1/2}}
  \sum_{y \in
W^{\sharp}/W} e^{-\pi i r  (g_{\Q}(\tilde y+\zeta,\tilde y +
\zeta))}.$$

\subsection{Volume  of lattices}\label{azed}  We shall use two simple properties of the
volume of lattices  defined in the Introduction.  Let $W$ be  a
lattice  of finite rank   and     $\langle \cdot, \cdot \rangle$
be an inner product in   $W_{\R}=W \otimes_{\Z} \R$. For a
sublattice   $W' \subset W$ of the same rank, the quotient  $W/W'$
is a finite abelian group and $\left| W/W' \right|=
 {\rm{vol}}(W') /{\rm{vol}}(W)   $.
 For the dual lattice $W^{\bullet}=\{ x \in W_{\R}: \langle x,W\rangle
\subset \Z \}$, we have  $
  {\rm{vol}}(W^{\bullet}) = ({\rm{vol}}(W))^{-1}$.

\subsection{Proof of (\ref{2eq:rec-ht})}
Note that we have two dual lattices $W^{\bullet}\subset W_\R$
(duality with respect to $\langle \cdot, \cdot \rangle$) and
$W^{\sharp} \subset W_\R$ (duality with respect to $g$). It
follows from the definitions that
 $W^{\bullet} = h(W^{\sharp})$.

It follows from the hypotheses $(\ref{eq:hypo})$ that if $r$ is
odd, then $g(x,x) \in 2\Z$ for all $x \in W$, so that $w=0$ is a
Wu class of $g$. If $r$ is even, we take $w$ to be an arbitrary
integral Wu class for $g$.
 Define $f$
and $v \in \Z$ as in Section  \ref{subsec:further} and set $\zeta
= h^{-1}(\psi)\in W_\R$. The hypotheses $(\ref{eq:hypo})$ imply
that $\zeta\in  \frac{1}{r}  W^{\sharp}$ as required in Section
\ref{subsec:further}. Therefore we have Formula
(\ref{eq:special1}). We can rewrite it  as
\begin{equation} \label{eq:special-almost-ht}
 \sum_{x \in W/rW}   e^{\frac{\pi i}{r}  (   \langle
\tilde x, h(\tilde x) \rangle + 2 r  \langle \tilde x, \psi
\rangle ) }
\end{equation}
$$
= e^{\frac{  \pi i}{4}\sigma(g) }{|G_g|^{-1/2}}\, r^{l/2}
 \, \sum_{y \in W^{\bullet}/h(W)}
e^{ -\pi i r \langle \tilde  y+\psi,h^{-1}(\tilde y+\psi) \rangle
}
 $$
where $\tilde x\in W$ is an arbitrary lift of $x$ and $\tilde y\in
W^{\bullet}$ is an arbitrary lift of $y$. By the assumptions,
$h(W^{\bullet}) \subset W^{\bullet}$. Observe that for $\tilde
y\in W^\bullet$, the expression $e^{-\pi i r \langle \tilde
y+\psi,h^{-1}(\tilde y+\psi) \rangle}$ depends only on $\tilde y
({\rm {mod}}\, h(W) + h(W^{\bullet}))$.  Indeed, for $z\in
h(W^{\bullet})\subset W^{\bullet}$,
$$r \langle z+\tilde y+\psi,h^{-1}(z+\tilde y+\psi) \rangle $$
$$= r \langle  \tilde y+\psi,h^{-1}( \tilde y+\psi) \rangle +2 r \langle \tilde y,h^{-1}(z )
\rangle  +2 r \langle  \psi,h^{-1}(z ) \rangle  +r \langle z
,h^{-1}(z ) \rangle .$$ By the assumptions (\ref{eq:hypo}), the
last three terms belong to $2\Z$. Therefore
$$e^{-\pi i r \langle z+\tilde y+\psi,h^{-1}(z+\tilde y+\psi) \rangle}= e^{-\pi i r \langle \tilde y+\psi,h^{-1}(\tilde y+\psi) \rangle}.$$
Set
$$s= \sum_{  y \in
W^\bullet/h(W)+h(W^{\bullet})} e^{-\pi i r \langle \tilde
y+\psi,h^{-1}(\tilde y+\psi) \rangle}.$$ Then
$$
\sum_{y \in W^{\bullet}/h(W)} e^{-\pi i r \langle \tilde
y+\psi,h^{-1}(\tilde y+\psi) \rangle}   = \Biggl| \frac{h(W) +
h(W^\bullet)}{h(W)} \Biggl| \, s  $$
$$= \Biggl| \frac{h(W^\bullet)}{h(W) \cap h(W^{\bullet})} \Biggr|
\, s  =  | h ( {W^{\bullet}}/{W \cap W^{\bullet}} )
 | \, s= |  {W^{\bullet}}/{W \cap W^{\bullet}}  | \, s.
$$
Similarly,
$$ \sum_{y \in W^{\bullet}/h(W^{\bullet})} e^{-\pi i r \langle
\tilde y+\psi,h^{-1}(\tilde y+\psi) \rangle} =  |  {W}/{W \cap
W^{\bullet}}  | \, s.$$ We deduce that
$$ \sum_{y \in W^{\bullet}/h(W)} e^{-\pi i r \langle
\tilde y+\psi,h^{-1}(\tilde y+\psi) \rangle}$$
$$=    |  {W^{\bullet}}/{W \cap W^{\bullet}} | \cdot
 |  {W}/{W \cap W^{\bullet}}  |^{-1} \sum_{y \in W^{\bullet}/h(W^\bullet )} e^{-\pi i r \langle
\tilde y+\psi,h^{-1}(\tilde y+\psi) \rangle}. $$ By the results of
Section \ref{azed},  $ {\rm{vol}}(W^{\bullet})= ( {\rm{vol}}(W
))^{-1}$,   $$ |  {W^{\bullet}}/{W \cap W^{\bullet}}|  \cdot |
{W}/{W \cap W^{\bullet}} |^{-1} =
\frac{{\rm{vol}}(W)}{{\rm{vol}}(W^{\bullet})}=({\rm{vol}} ( W ))^2
$$ and  $$ |G_{g}|= |W^{\sharp}/W| =  |W^\bullet/h(W)| =
{{\rm{vol}} (h(W))} / {{\rm{vol}} (W^\bullet)} =  |\det(h)| \,
({\rm{vol}} ( W ))^2.$$
  Substituting these formulas in (\ref{eq:special-almost-ht}),
we obtain Formula   (\ref{2eq:rec-ht}).  \hfill $\blacksquare$\\

\noindent{\bf{Remark}}. The proof above goes mutatis mutandis for
any  nondegenerate symmetric bilinear  pairing $\langle \cdot,
\cdot \rangle$ on $W_\R$. (The arguments invoked in
\cite{HansenTakata} and \cite{Jeffrey} require the pairing to be
positive definite.)

\bibliographystyle{amsalpha}

\vskip 0.5cm

\noindent{\emph{Adresses of authors:}}\\

\small

\noindent F.D., Institut de Math\'ematiques, Universit\'e Paul
Sabatier -- Toulouse III, 118, route de Narbonne, 31062 Toulouse
cedex 4, France. Email:
deloup@picard.ups-tlse.fr \\

\noindent V.T., Institut de Recherche en Math\'ematiques
Avanc\'ees -- UMR 7501 CNRS , Universit\'e Louis Pasteur --
Strasbourg I, 7 rue Ren\'e Descartes, 67084 Strasbourg cedex,
France. Email: turaev@math.u-strasbg.fr

\end{document}